\input amstex
\input epsf
\input amsppt.sty
\NoBlackBoxes
\magnification=1200
\magnification=1200 \vsize=7.50in\parindent 20 pt
\NoBlackBoxes
\nopagenumbers
    \define \Dl{\Delta} \define \dl{\delta}
          \define
\vp{\varphi}

\define \ri{\rightarrow}

 \define \ep{\endproclaim}

  \define \bk{\bigskip}

\define \1{^{-1}} \define \2{^{-2}}

   \define \BC{\Bbb C}  
  \define \BP{\Bbb P}

\define\Center{\operatorname{Center}}

\baselineskip 20pt

\parindent 20 pt
\NoBlackBoxes


\topmatter

\title BMT Invariants  of Surfaces  and 4-Manifolds\endtitle
\author Mina Teicher\endauthor
\address Department of Mathematics and Computer Science, Bar-Ilan
University, 52900 Ramat-Gan, Israel\endaddress
\email teicher\@macs.biu.ac.il\endemail
\thanks This research is partially supported by the Emmy Noether Research
Institute for
Mathematics and by grant no. 8007199 of the Excellency Center for
Group-Theoretic
Methods in Algebraic Varieties of the Israel National Academy of
Sciences and by EAGER, the European Network in Algebraic
Geometry Education and Research.\endthanks
\endtopmatter

\document
\baselineskip 20pt
\subheading{Introduction}

In this paper we present the Braid Monodromy Type (BMT) of curves and
surfaces; past, present and future. The BMT is an invariant that can
distinguish between non-isotopic curves; between different families
of surfaces of
general type;  between connected components of moduli space of
surfaces and between non symplectmorphic
4-manifolds.
 BMT is a  finer invariant than the Sieberg-Witten invariants.  Consider $X_1$,\ $X_2$ surfaces of general type with the same $c_1^2,c_2$\
$(\pi_1=1).$
It is known that
$$\alignat2& &&X_1\ \text{is a deformation of}\ X_2\\
&\Rightarrow &&X_1\ \text{is diffeomorphic}\ X_2\\
&\Rightarrow &&X_1\ \text{is homeomorphic to}\ X_2\\
& \Rightarrow &&X_1\ \text{is homotopic to}\ X_2\endalignat$$

What about the reverse directions?
Are there invariants distinguishing between these equivalence classes?
The new invariant, proposed here, is located between the first and the second
arrow. In this paper we shall introduce the new invariant, state the current
results and pose an
open question.

We start by defining in Section 1 the Braid Monodromy Type (BMT) of a curve.
If $S$ is a curve of degree $m=\deg S,$
and $\rho: \pi_1(\Bbb C-N)\to B_m[\BC_u,\BC_u\cap S]=B_m$ is the braid
monodromy (BM) (see Section 1) of $S,$
then
it is known (see for example \cite{MoTe3}) that if $\{\dl_i\}$ is a $g$-base of
$\pi_1(\BC-N),$\ $\Dl^2$ is the generator of $\Center(B_m)$, then (Artin):
$\Dl^2=\Pi\rho(\dl_i).$ Such a factorization is called a braid monodromy
factorization of $\Dl^2$ related to $S.$
We define a (Hurwitz) equivalence relation on (positive) factorizations of
$\Dl^2.$
We prove (see, \cite{MoTe3}) that all factorizations of $\Dl^2$ induced
from a curve $S$ are
equivalent and occupy a full equivalence class.

A BMT of a curve is the equivalence class of factorizations.
Together with V.~Kulikov we proved in 1998 that a BMT of $S$ determines
isotopy type of $S.$
We want to use the BMT of the branch curve of a surface   as an invariant of the surface.
First, we have to answer the following question:\ Are surfaces determined
by their  branch
curves?
$$\alignat5
&(\ X && , \quad f\ )\quad&&\qquad \   \binom{?}{\Leftrightarrow}\quad
&&\qquad\quad S\\ &\quad \downarrow  &&\quad\downarrow\quad&&\qquad\quad
\downarrow\quad&&\qquad \quad\downarrow\\
&\text{surface}&&\quad\text{generic projection to}\
\BC\BP^2\quad &&\text{Chisini-Kulikov 97}&&\quad\text{branch
curve}\endalignat$$

This question was known  as the Chisini Conjecture and was proved by
V.~Kulikov in 1997 (for
$\deg S\gg0) $ (see \cite{Ku2}).
The BMT of branch curves is indeed  an invariant of surfaces.
Moreover, this invariant is very powerful.

\medskip

The paper is divided as follows:

{\S1.  Braid Monodromy Type (BMT) of Curves}

\S2. Braid Monodromy Type (BMT)  of Surfaces

\S3. How to Compute the Braid Monodromy Type (BMT) of a Surface

{\S4. Future plans}

\bk

\subheading{\S1. Braid Monodromy Type (BMT) of Curves}

We shall start by defining the braid monodromy.
We shall derive a discrete invariant of $S$ from the braid monodromy of $S$
and its
factorization.

\definition{Definition} {\it The Braid Monodromy}

Let $S$ be a curve, $S\subseteq \BC^2,$
\ $\Pi(x,y)=x,$\ $\Pi$ is the projection on the $x$-axis,
\
$(\deg\Pi=m),$\ $\Pi$ rstricted on $S.$
Let $N=\{x\in\BC^1\bigm| \#\Pi\1(x)\lneqq m\}.$
Take $u\notin N,$ s.t. $x\ll u,\ \forall\ x\in N$ and $\#\BC_u^1\cap S=m.$

 {\it The braid monodromy with respect to} $S,\Pi,u.$ is a naturally
defined homomorphism:
$$\pi_1(\BC^1-N,u)\overset\vp\to\ri B_m[\BC_u^1,\BC_u^1\cap S]$$
 \enddefinition

\remark{Remark} The classical monodromy of the cover $S\rightarrow (x$-axis factors
through the braid monodromy
$$\alignat 3
&\pi_1\quad && \overset\vp\to\ri &&\quad B_m\\
& \quad&&\searrow &&\quad \downarrow\\
& && &&\quad S_m\endalignat$$\endremark

We recall Artin's theorem concerning a presentation of 
the generator of center of the braid group (which is known also as the Dehn-twist of the braid group )
 as a product of braids which are the images of the elements of a geometric base under the braid monodromy.

\proclaim{Theorem} \rom{(Artin)}

$\vp$ is a braid monodromy of $S,$\ $\vp:\pi_1\to B_m.$  Let $\{\dl_i\}$ be
a geometric
(free) base ($g$-base) of
$\pi_1.$ (see figure 1).
 $\Dl^2$ is the generator of $\Center(B_m).$
Then:
$$\Dl^2=\Pi \vp(\dl_i)$$.\ep

\centerline{
\epsfysize=4cm
\epsfbox{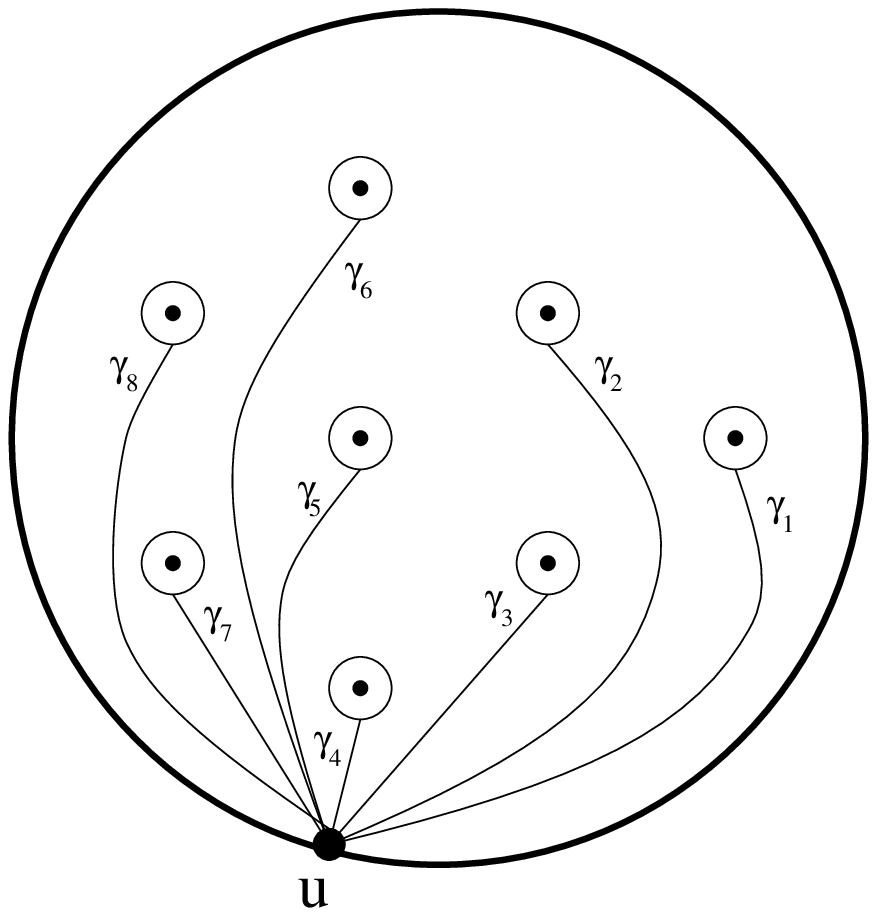}
}

\medskip

Following the theorem, we  define
\definition{Definition} {\it Braid monodromy factorization (BMF) related to
a curve $S$ and a geometric-base
$\{\dl_i\}$}

A presentation of $\Dl^2$ as a product of the form $\Dl^2=\Pi\rho(\dl_i)$ for
$\dl_i$ a $g$-base.\enddefinition

Clearly, there are many BMF related to a curve; each geometric base will induce another BMF of the curve.
In fact, a BMF is a special case of positive factorization of $\Dl^2:$

\definition{Definition} {\it Positive factorization  (PF) of $\Dl^2$}

A presentation of $\Dl^2$ as a product of  (by a positive braid we mean all braids that can be presnted as conjugation of 
frame elements).\enddefinition

For example, $\Dl_3^2=X_1X_2X_1X_2X_1X_2$ is a PF.  Artin proved that a BMF
is a PF (clearly,
he used different terminology). In 
order to discuss an equivalent relation on BMF as we do next , we first have to consider the set of PF's in general .
Later it will be clear that the equivalent relation we define on the PF's is closed on BMF's.
We shall define an equivalence relation on PF's using Hurwitz moves.

\definition{Definition} {\it Hurwitz move (HM)}

$a_i,b_i\in G$ is a group. \ $(a_1,\dots,a_n)$ is obtained from
$(b_1,\dots,b_n)$ by an HM
if $\exists k$ s.t. $$\gather a_i=b_i,\ i\ne k,k+1\\
a_k=b_kb_{k+1}b_k\1\\
a_{k+1}=b_k\endgather$$\enddefinition

\definition{Definition} {\it Equivalent PF's}

$\Dl^2=\Pi a_i$ is equivalent to $\Dl^2=\Pi b_i$ if $(a_1,\dots,a_n)$ is
obtained from
$(b_1,\dots,b_n)$ by a finite number of HM's.\enddefinition

\proclaim{Theorem} (\cite{Mote3})
\roster\item"(a)" Two BMF's related to $S$ (for different $g$=bases 
$\{\dl_i\}$) are equivalent.
\item"(b)" If a PF is equivalent to a BMF, then the PF is also a BMF (with
another
$\{\dl_i\}$ and the same $S).$\endroster\ep

\proclaim{Corollary} The set of BMF's related to a curve $S$ and different
$g$-bases
(different $\{\dl_i\}$'s) occupy a full \rom{equivalence class} of PF's.\ep

\remark{Remark} By the previous corollary, if there exist two BMF's related
to $S$ which
are equivalent, then  any two BMF's are equivalent.\endremark

Using the above theorem, we define
\definition{Definition} {\it Braid monodromy type of curves (BMT)}

Two curves $S_1$ and $S_2$ are of the same BMT if they have related BMF's
that are
equivalent.\enddefinition

\medskip
It can be deducted from 
Moishezon work (1992) 
on the counterexample to Chisini theoreme that for any BMT 
with algebraic factors there is a semi-algebraic curve having this BMT. 
This does not determins uniqueness, which was donein 1998.

In 1998, we proved that given 2 semi-algebraic curves $S_2$ and $S_1$ ,
if $S_1\overset{\text{BMT}}\to\cong S_2$ then $S_1$ is isotopic to $S_2$. We 
proved:

\proclaim{Theorem 1} \rom{(Kulikov, Teicher)\ (AG9905149)}

If $S_1$ and $S_2$ are two semi-algebraic cuspidal curves 
which have the same BMT, then 
there is a symplectic isotopy $ F_t: CP^2\to CP^2$ s.t. 
(1) $F_t(S_1)$ is a Hurwitz curve for all $t$; 
(2) $F_1(S_1)=S_2.$ .\ep

\bk

\subheading{\S2. Braid Monodromy Type (BMT) of Surfaces}

We shall derive an invariant of a surface $X$ from the braid monodromy of
its branch curve
$S$ related to a generic projection $X\to\BC\BP^2$\ $(S\subseteq
\BC\BP^2)$.

\medskip

\definition{Definition} {\it Braid monodromy type of surfaces (BMT)}

The BMT of a surface of general type is the BMT of the branch curves of a
generic projection of the surface imbedded in a projective space by means
of the
multiple canonical class.\enddefinition

This definition can be extended to any 4-manifold for which one can
construct a ``good''
generic projection to $\BC\BP^2$ (see, for example, a construction of a
projection of a
symplectic 4-manifold  in Aroux 
\cite{A}).

\newpage

We proved in 1998:
\proclaim{Theorem 2} \rom{(Kulikov, Teicher) \ (AG9807153) (AG9905129)}

The Braid Monodromy Type (BMT) of a projective surface $X$, which
is a  generic covering 
s.t. Chisini's conjecture holds for its branch curve,  determines
the diffeomorphism type
of
$X.$\ep

This  theorem is determining the position of the BMT invariant of surfaces relative the 
known invariants and might give the final invariant that 
will determine the Diff-Def question: How 
to distinguish among diffeormophic surfaces which are not deformation of each other? 

It remains an open question whether the inverse is correct or whether the
Braid Monodromy Type
(BMT) of
$S$ determines the deformation type of $X.$
Clearly,
the BMT of a surface determines the fundamental group (f.g.) of the
complement of the branch
curve.
In the following diagram we illustrate the relations between different
invariants of
surfaces.
$$\alignat 4
& &&\binom{?}{\Leftarrow}&&  &&   \binom{?}{\Leftarrow} \quad\ \
\not\Leftarrow\\
&\text{Def} && \ \Rightarrow \  &&\text{BMT}\ &&
\ \Rightarrow \ \text{Diff}\
\Rightarrow\ \text{Hom}\\
& && &&\ \Downarrow \\
& && && \ \text{f.g.}\endalignat$$

The proof of Theorem 2 uses Theorem 1 and the Chisini Conjecture from 1950
(= Kulikov Theorem
1998). If we have two surfaces and two generic projections and 2 branch curves
$$\alignat 2
& && X_i\\
& && \downarrow\\
&S_i\subset\ && \BC\BP^2\endalignat$$
  (with $\deg S_i>e$), then Kulikov proved that $S_1\cong S_2\Rightarrow
X_1\cong X_2.$
  Kulikov also proved, (see \cite{Ku1}), that
$$\align\left\{\aligned  \qquad X\ \\
  \downarrow\ \\
S\subset \BC\BP^2\endaligned\right\}
\overset{1-1}\to\cong\{\Pi_1(S^c)\twoheadrightarrow\text{symmetric
group}\},\qquad\qquad\qquad
\endalign$$
and we use this, too.  
\bk
\subheading{\S3. How to Compute the BMT of a Surface}

To compute the BMF and the BMT of branch curves, we use degenerations of
the surfaces to
union of planes.

Let $X$ be a surface of general type.  We degenerate it  into $X_0.$
$$\alignat 3
& &&X\ \rightsquigarrow\ && X_0= \text{union of planes}\\
& &&\downarrow &&\downarrow\\
&S\subset\ &&\BC\BP^2 && \BC\BP^2\supset S_0,\qquad S_0\ \text{is union of
lines}\endalignat$$

It is difficult to derive $S$ from $S_0$ (we know that $\deg S=2\deg S_0$).
Instead,
we derive the BMF of $S$ from the BMF of $S_0.$
This is called the {\it Regeneration Process}.
It is described in detail in \cite{MoTe4}, \cite{RoTe}, \cite{AmTe}.

The regeneration process of BMF is divided into 4 main steps as follows:

\roster\item"$\bullet$" {\it Computing the braid monodromy of line
arrangements}
(\cite{MoTe3})
\item"$\bullet$" {\it Microscopic techniques}\newline Regeneration of
arrangements with one
singular point (when 2 lines or 3 lines or 4 lines (etc.) intersect at one
point. The
computation for $n$ lines meeting at one point is based on the  computation
of the
regeneration of
$n-1$ lines which is based on  the previous one.
\item"$\bullet$ $\bullet$ $\bullet$" {\it Global arguments} By comparing the sum of the degrees of the all factors 
arrising 
from the degeneration process  to the known degree of the generator of the center 
(which is n(n-1)) in order to determine how many missing factors we oughtt o look for and of what degrees.
\item"$\bullet$ $\bullet$" {\it Counting extra branch points}\newline
Looking for branch points close to infinity.
(singularities not arising from a regeneration of a singular point). \endroster
\bk
\roster\item"$\bullet$" The {\it microscopic techniques} for branch curves
include 3 basic
regeneration rules:\newline
I.\quad Regeneration of a branch point  (one branch point is replaced by 2
branch\linebreak
{}\quad points).
\newline II.\ \ Regeneration of a node (one node is replaced by 4 nodes).
\newline
III. Regeneration of a tangent point (a tangency point is replaced
by\linebreak {}\quad 3
cusps).\endroster

\bk\bk\bk

(See \cite{MoTe4}, \S3, Lemmas 3.1, 3.2 and 3.3.)

\subheading{\S4. Future Plans}

As par

\roster\item Finding an algorithm that determines when two BMF's are Hurwitz 
equivalent  i.e., if one can be derived from the other by a finite number 
of Hurwitz moves. 
This work is being carried out in the framework of sub-projects
concerning:\ the word
 problem in the braid group 
and its complexity ,
Hurwitz equivalence of words made of frame elements 
and computerized algorithms; Completing this project is complementary to:

\item Looking for any easy computable discrete invariant 
of BMT of 
a curve which is easier to compute than the
fundamental group of the complement of the  curve ;

\item Finding an algorithm that determines if a certain PF is a BMF related
to a curve 
and in particular to a branch curve; 

This is 
essential in order 
to derive all equivalence classes of BMT of branch curves which will result
in  determining deformation classes of surfaces. The next 2 items deal with
surfaces and the last one with 4-manifolds.

\item Finding degeneration techniques into union of planes of
different  surfaces ;
This will allow us to compute 
braid monodromy of branch curves via regeneration techniques as in 
the next item:
\item Computing regeneration of an intersection of $n \ (\gg 1)$ lines;

\item Adapting the BMF to symplectic 4-manifolds (after Aroux);
\item Computing related fundamental groups; Galois groups, groups of line 
arrangemnts, complements of branch curves of  triple covers .\endroster

\enddocument